\documentclass[11pt,twoside,a4paper,notitlepage]{article}

\usepackage{amsmath}

\usepackage{amssymb}
\usepackage{amsfonts}
\usepackage[a4paper,margin=2cm,vmargin={1cm,2cm},includeheadfoot]{geometry}
\usepackage[T1]{fontenc}
\usepackage{amsthm}
\usepackage{enumerate}
\usepackage[scaled=0.92]{helvet}
\usepackage{srcltx}
\usepackage{graphicx}
\usepackage{tikz}
\usepackage{booktabs}
\usepackage{subfigure}

\usepackage[labelfont=bf,margin=1cm,justification=justified,labelsep=period]{caption}

\theoremstyle{plain}
\newtheorem{theorem}{Theorem}[section]

\newtheorem{lemma}[theorem]{Lemma}

\newtheorem{property}[theorem]{Property}
{\theoremstyle{definition}
\newtheorem{remark}[theorem]{Remark}
\newtheorem{notations}[theorem]{Notations}
\newtheorem{remarks}[theorem]{Remarks}}

\newcommand{\R}{{\mathbb {R}}}

\newcommand{\bs}{\boldsymbol}   
\newcommand{\wt}{\widetilde}    
\newcommand{\wh}{\widehat}      

\newcommand{\buil}[3]{\mathrel{\mathop{\kern0pt#1}\limits_{#2}^{#3}}}

\definecolor{gris}{gray}{0.5}

\begin{document}

\begin{center}
{\LARGE\bf Maxima of
 Weibull--like\\ [10 pt] distributions and the Lambert W function} \\ [0.5 cm]

{\sc\large Armengol Gasull \footnote{Department of Mathematics, Universitat Aut\`{o}noma de
Barcelona, {\it e--mails:} gasull@mat.uab.cat, utzet@mat.uab.cat},
 Jos\'e A. L\'opez--Salcedo \footnote{Department of Telecommunications
  and Systems Engineering,
Engineering School, Universitat Aut\`{o}noma de
Barcelona, {\it e--mail:} jose.salcedo@uab.cat }, Frederic Utzet $^1$} \\ [0.5 cm]


\end{center}

\noindent{\bf Abstract.} The Weibull--like distributions form a large class of probability
distributions
 that belong  to the domain of attraction for the
maxima of the Gumbel law. Besides the Weibull distribution, it includes
 important distributions as  the Gamma laws and, in particular,
the $\chi^2$ distributions. In order to have explicit expressions of
the  norming  constants for the maxima
 it is necessary to solve asymptotically  a nonlinear
equation; however, for some members of that family,   numerical and  simulation studies show that the
   constants that are usually
suggested  are  inaccurate for moderate or even large sample sizes.
 In this paper we propose  other norming constants computed
with the asymptotics of the Lambert W function that significantly    improve   the
accuracy of the approximation to the Gumbel law.
These results are applied  to the computation of the constants for   the maxima of
Gamma random variables that appear in some applied problems.



\section{Introduction}
The origin of this paper was when the authors tried to use Extreme Value Theory to the
maxima of $\chi^2$ random variables in an applied problem of signal processing
(see Turunen \cite{Tur07}).
 In particular, the goal was to accurately characterize the detection performance
  of a Global Navigation Satellite System (GNSS) receiver,
   whose main task is to provide positioning information by processing
   the signals emitted from Earth--orbiting satellites (see Seco--Granados {\it et al.}
   \cite{Seco12}).
  To do so, a GNSS receiver must first detect the presence of visible satellites,
  which is done by analyzing the signal that impinges onto the GNSS receiver antenna.
  This analysis requires a bi--dimensional search in order to determine the carrier
   frequency and the time--delay for each of the satellite signals of interest,
   in a similar manner to what occurs when tuning a radio into a specific radio
    station and a specific program, respectively. For each of the tentative carrier
     frequency and time-delay values, the GNSS receiver measures the received signal
    power, which can be modeled as a $\chi^2(m)$ random variable, and stores the
     resulting power measurement into a specific cell of a time-frequency matrix
      (see   Seco--Granados {\it et al.} \cite[Eq. (9)]{Seco12}). When all tentative values
       have been tested, the GNSS receiver takes the maximum of all
       the entries within the time-frequency matrix. This
        leads to the so-called ``parallel acquisition'' approach,
        and the resulting maximum value is then compared to a threshold
         in order to determine whether the satellite being analyzed was
          actually present or not (see Seco--Granados {\it et al.}
 \cite[Eq. (11)]{Seco12}).

In practice, this parallel acquisition approach typically
 involves computing the maximum of 10 to $10^6$ $\chi^2(m)$ random variables, for
 $m$ between 10 and 20,
  depending on whether assistance information is provided or not to the GNSS
   receiver. In this context, it is interesting to note that $\chi^2$
    random variables are in the domain of attraction for maxima of the Gumbel law,
     however, the norming constants that are usually proposed give,
      for such  sample sizes, inaccurate results (see Subsection \ref{sec:constant-chi}
      ).  Then  we realized that the computations needed to
obtain these constants were related with the Lambert W function
 (see Corless {\it et al.}  \cite{CorGonHarJefKnu86}), and that the asymptotic expansion  for that
 function helps very much to improve the norming constants. As Resnick \cite[Page 67]{Res87}
 points out,  {\it Computing normalizing constants can be a brutal business, and
 any techniques which aid in this are welcome indeed}.
 The purpose of this paper
 is to show how the asymptotics of the Lambert W function and its generalizations can
 be used in this problem.

 We show that the {\it centering} constant for the maxima of $n$ i.i.d.  random
 variables with  distribution called
 Weibull-like (see Section \ref{sec:intro}) can be expressed
 in terms of the secondary branch of a real Lambert W function; the asymptotic expansion of
 that function is well known, and the centering  constant that is deduced using standard methods
 of asymptotic analysis corresponds to the two dominant terms of that expansion,
 loosely speaking,  of the form $C_1 \log n+C_2 \log\log n$.
 However, for typical sample sizes  the results
 are quite inaccurate, and we propose to add one more term of that asymptotic
 expansion, basically of the form $\log\log n/\log n$; this term goes to zero when
 $n\to\infty$, but so slowly
 that cannot be neglected. We pay special attention to the maxima of Gamma laws.
 In that case, we need a double enhancement of the standard technique: on the one
 hand, the usual distribution tail equivalent to a Gamma distribution needs to be improved;
 on the other hand, using the asymptotic expansion of a generalization of
 Lambert W function,  we add an additional term, that, as before,  goes to zero when the sample size
 increases, but also helps very much to get accurate approximations.

The contents of the paper are the following. In Section 2 we introduce the class
of generalized  Weibull distributions and its tail-equivalent distributions, called
Weibull-like distributions; we also recall some essential facts  about Extreme
Value Theory. In Section 3 we comment the main results of the paper. In Section
 4 we study a particular simple case
of a generalized Weibull  distribution and we describe the problem of the inaccuracy
of the norming constants; further we introduce the
Lambert W function and its asymptotics.  In Section 5 we study the velocity
of convergence to the maxima and show the importance of the election of
the norming constants. In Section 6 we apply these results to the maxima
of Gamma laws.  Some technical matters are placed in the Appendix.

\section{Weibull--like distributions and their maxima}
\label{sec:intro}
To introduce  notation and to describe the context of the paper we recall
a few basic facts from Extreme Value Theory.
Let  $X_1,\dots,X_n$ be i.i.d. random variables with common cumulative  distribution function $F$,
and denote by $M_n$ its maximum,
$$M_n=\max\{X_1,\dots,X_n\}.$$
It is said that $F$ is in the domain of attraction for maxima of the Gumbel law
if
 there are sequences of real numbers $\{a_n,\ n\ge 1\}$ and $\{b_n,\ n\ge 1\}$
(the {\it norming --or normalizing-- constants})
with $a_n>0$ such that
\begin{equation}
\label{conv0}
\lim_n \frac{1}{a_n}\,(M_n-b_n)= H, \ \text{in distribution},
\end{equation}
where $H$ is a Gumbel random variable, with  distribution function
\begin{equation}
\label{gumbel-dist}
\Lambda(x)=\exp\{-e^{-x}\},\ x\in \R.
\end{equation}
The norming constants can be taken
 (see,
for example,  Resnick \cite[Proposition 1.11]{Res87})
\begin{equation}
\label{exacta-weibull-pre}
b_n=F^{-1}(1-n^{-1})
\end{equation}
and
\begin{equation}
\label{exacta-weibull-c-pre}
a_n=A(b_n)
\end{equation}
where $A(x)$ is an auxiliary function of the distribution function $F$.
Auxiliary
functions are not unique though they are asymptotically  equal.
However, under certain conditions (in particular,  $F$ should have   density, denoted by $f$, for $x>x_0$, for
some $x_0$)  an  auxiliary
function  is (see  again Resnick \cite[Proposition 1.11]{Res87})
\begin{equation}
\label{aux-canon} A(x)=\frac{1-F(x)}{f(x)}.
\end{equation}
We should remark that from the standard proof of the convergence (\ref{conv0}) it is not
deduced that these constants  produce more accurate results than
other constants computed with other auxiliary functions or other
ways.

 In order to obtain explicit expressions of $b_n$ and $a_n$
 the following two results are  used: the first one can be called  {\it simplification
by tail equivalence} (see Resnick \cite[Proposition 1.19]{Res87}):

\begin{property}
\label{propietat1}
 Let $F$ be a distribution function in the domain of attraction
of a Gumbel law, and let $G$ be another distribution function  right tail equivalent
 to $F$:
$$\lim_{x\to \infty}\frac{1-G(x)}{1-F(x)}=1.$$
Then $G$ is also in the domain of attraction of the Gumbel law and
  the norming constants  of  $F$ and $G$ can be taken equal.
\end{property}

The second property is just a property of the convergence in law applied to
our context.

\begin{property}
\label{propietat2}
Let $F$ be a distribution function that belongs to the
domain of attraction for maxima of the Gumbel law, with
norming constants
  $\{a_n,\ n\ge 1\}$ and $\{b_n,\ n\ge 1\}$.
 If the  sequences
$\{a_n',\, n \ge 1\}$ and $\{b_n',\, n \ge 1\}$ satisfy
$$
\lim _n \dfrac{a_n}{a'_n}=1\quad \text{and}\quad
 \lim_n\dfrac{b_n-b'_n}{a_n}=0,$$
 then
$$\lim_{n}\frac{1}{a'_n}\big(M_n-b'_n\big)=H\ \text{in distribution,}$$
that is, the sequences
$\{a_n',\, n \ge 1\}$ and $\{b_n',\, n \ge 1\}$ are also norming constants for $F$.
\end{property}

We will study a rich  family of distribution functions. To begin with,  in this paper, a probability distribution function $F$ such that for some $x_0$
has the form
\begin{equation}
\label{weibull-type}
F(x)=1-Kx^\alpha \exp\{-Cx^\tau\},\ x\ge x_0,
\end{equation}
where $K,\,C,\,\tau>0,$ and $\alpha\in\R$, will be called a
{\it generalized Weibull distribution}; the standard Weibull  law $W(\lambda,\nu)$,
where $\lambda>0$ is the scale parameter and $\nu>0$
 the shape parameter,  is the case $\alpha=0,\, \tau=\nu,\ C=1/\lambda^\nu, \, K=1$
and $x_0=0$; in particular an exponential law of parameter $\lambda>0$ has
$\alpha=0,\, \tau=1$ and $C=\lambda$, and  a $\chi^2(2)$ law has $\alpha=0, \, \tau=1,\, C=1/2$.
It is easy to check that generalized Weibull distributions belong to the domain of atraction for maxima to the Gumbel law.

In agreement  with Embrechts {\it et al.} \cite[page 155]{EmbKluMik97},
a probability distribution function  right tail equivalent
to a generalized Weibull distribution
is said to be a {\it Weibull--like distribution}.
A main example is
 the Gamma law $G(\nu,\theta)$ (with $\alpha=\nu-1,\ \tau=1$), and, in particular,
 a $\chi^2(m)$ law, see Section \ref{chi2}.
The normal law also is Weibull-like with  $\alpha=-1$ and $\tau=2$; however, this case
has specials properties:  on the one hand,
there is the remarkable result of Hall  \cite{Hall79} where  he proves that
for some norming constants $a^*_n$ and $b^*_n$,
$$\frac{C_1}{\log n}<\sup_{x\in \R}\vert \Phi^n (a_n^*x+b_n^*)-\Lambda(x)\vert
 < \frac{C_2}{\log n},$$
($\Phi$ is the cumulative distribution function of the standard normal law) where $C_2$ can be taken equal to 3, and that the rate of convergence cannot be improved by choosing a different sequence
of norming constants; on the other hand, the fact that $\alpha<0$ introduces important
changes in our approach; for some improvements on the norming constants for the normal case, see
Gasull {\it et al.} {\cite{GasullJolisUtzet}.
 Given that we are mainly interested in the Gamma law
 we will assume from now on that $\alpha>0$ and $\tau\ge 1$.

Consider $X_1,\dots,X_n$ be i.i.d. random variables with  Weibull-like distribution function
 $G$, right tail equivalent to a generalized Weibull distribution function $F$ of the
 form (\ref{weibull-type}).
Thanks to Property \ref{propietat1},
the norming constants can be taken
\begin{equation}
\label{exacta-weibull}
\begin{array}{l}
b_n=F^{-1}(1-n^{-1})\\
a_n=\dfrac{1}{C\tau b_n^{\tau-1}-{\alpha}/{b_n}}
\end{array}
\end{equation}
where  for the expression for $a_n$ we have used the  auxiliary function (\ref{aux-canon}) associated
to $F$. From that expressions, by  using Property
\ref{propietat1} and asymptotic analysis, it is possible to find explicit expressions
of the norming constants, and usually  are suggested,
see, for example, Embrechts {\it et al.} \cite[page 155]{EmbKluMik97},
\begin{equation}
\label{aprox-weibull}
\begin{array}{l}
b'_n=\big(C^{-1}\log n\big)^{1/\tau}+\dfrac{1}{\tau}\big(C^{-1}\log n\big)^{1/\tau-1}
\bigg(\dfrac{\alpha}{C\tau}\log\big(C^{-1}\log n\big)+\dfrac{1}{C}\log K\bigg)\\
a'_n=(C\tau)^{-1}\big(C^{-1}\log n\big)^{1/\tau-1}
\end{array}
\end{equation}
These will be called  the {\it standard  constants}.

\section{Main results}

Our purpose is to show that, for moderate or even quite large sample sizes,
 the election of the norming constants plays a major
role in the velocity of convergence of (\ref{conv0}),
and that it is possible to choose  constants that produce more accurate results
than   the standard ones. Our main finding is that  for a generalized Weibull
distribution  (\ref{weibull-type}), rather than the standard constants
(\ref{aprox-weibull}), it is more convenient to use other ones based on
 the asymptotic expansion of the Lambert W function
(see Subsection \ref{sub:lambert}). They are:

\bigskip

\noindent $\bullet$ If $\alpha> \tau$,

\begin{equation}
\label{bsecgeneral}
b''_n=\Big(\dfrac{\alpha}{C\tau}\Big)^{1/\tau}\bigg(-M_1+M_2-\dfrac{M_2}{M_1}\bigg)^{1/\tau},
\end{equation}
and
\begin{equation}
\label{asecgeneral}
a''_n=\dfrac{1}{C\tau (b''_n)^{\tau-1}-{\alpha}/{b''_n}},
\end{equation}
where
$$M_1=\log\Big(\frac{C\tau}{\alpha (Kn)^{\tau/\alpha}}\Big)\quad \text{and}
\quad M_2=\log(-M_1).$$

\bigskip

\noindent $\bullet$ If $\alpha\le \tau$,
\begin{equation}
\label{bsecalphapetit}
b_n''=\frac{1}{C^{1/\tau}}\Big(N_1+\frac{\alpha}{\tau}\, N_2+
\frac{\alpha^2}{\tau^2}\,\frac{N_2}{N_1}\Big)^{1/\tau},
\end{equation}
where
$$N_1=\log\Big(\frac{Kn}{C^{\alpha/\tau}}\Big)\quad \text{and}\quad
N_2=\log(N_1),$$
and $a''_n$ the same as in (\ref{asecgeneral}).

Although $b'_n$ of (\ref{aprox-weibull}) and these  $b''_n$ look like quite different,
$b'_n$ coincides, except some constants, with the first two terms of $b''_n$
(\ref{bsecgeneral})
(related with $M_1$ and $M_2$) or (\ref{bsecalphapetit}),
and the remaining part of $b_n''$ is  a sequence that converges
to zero, but so slowly that, jointly with the constants,
 it is important in typical  sample sizes.
See below the case of Gamma random variables, where the difference between the constants
is more evident.
To illustrate this point, in Section \ref{sec-casfacil} we study the simplest case of
a non--trivial generalized Weibull  distribution:
\begin{equation}
\label{dist-facil-main}
F(x)=\begin{cases}
1-e\, x\,e^{-x},& \text{if $x\ge 1,$}\\
0, & \text{if $x<1$,}
\end{cases}
\end{equation}
where  $K=e$ and $ C=\tau=\alpha=1$. The standard constants (\ref{aprox-weibull})  are
\begin{equation*}
\begin{array}{l}
b'_n=\log n+\log \log  n+1,\\
a'_n=1.
\end{array}
\end{equation*}
The proposed constants are
\begin{equation*}
\begin{array}{l}
b''_n=\log  n+\log  (\log  n+1)+1 +\dfrac{\log ( \log  n+1)}{\log  n+1},\\
a''_n=\dfrac{b''_n}{b_n''-1}.
\end{array}
\end{equation*}

The difference
between $b'_n$ and $b''_n$ is essentially $\log\log n/\log n$, that goes to zero,
 but so slowly
that cannot be neglected  in typical cases. This is illustrate in Figure \ref{dosHistogrames}.
Random variables from the distribution  (\ref{dist-facil-main}) are easily simulated
(see Subsection \ref{simulacio}). In that figure there are the histogram
of a simulation of $10^4$ maxima of $n=100$ random variables normalized by using
the standard norming constants $a'_n$ and $b'_n$, and by using the proposed
constants $a''_n$ and $b''_n$,

\begin{figure}[htb]
\centering
\subfigure[With standard norming constants $a'_n$ and $b'_n$]{%
\includegraphics[width=0.45\linewidth]{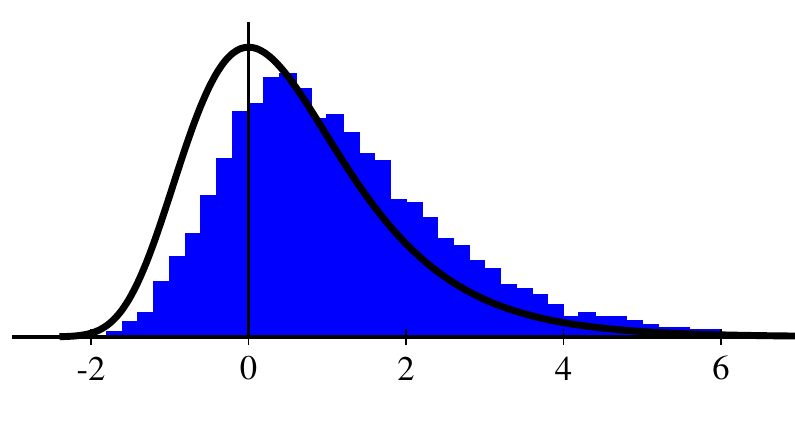}
       }%
\hspace{0.01\linewidth}
\subfigure[With proposed norming constants $a''_n$ and $b''_n$]{%
\includegraphics[width=0.45\linewidth]{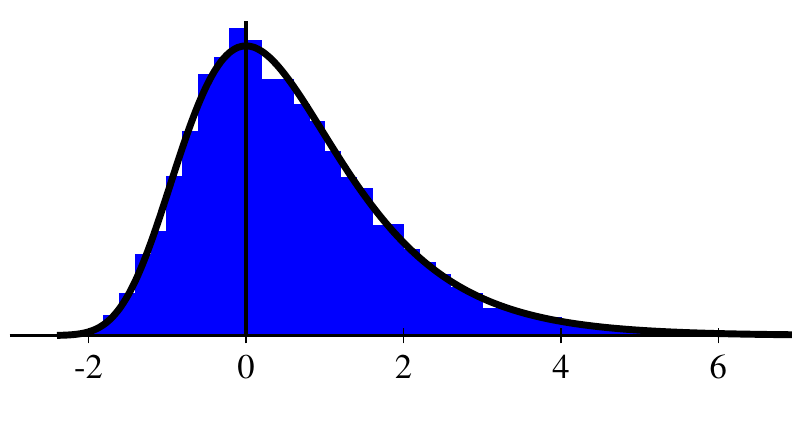}
}%
\caption{ Solid line:
Gumbel density. Histograms of a simulation of  $10^4$    maxima with $n=100$, of a generalized Weibull distribution
of parameteres $C=\tau=\alpha=1$, with two different sets of norming constants}
\label{dosHistogrames}
\end{figure}

\bigskip

  As we commented, our main interest is in  Gamma laws
$G(\nu,\theta)$ (with $\nu>0$ and $\tau\ge 1$); the standard constants (\ref{aprox-weibull})  are
(see Embrechts {\it et al.} \cite[page 156]{EmbKluMik97})
\begin{equation}
\label{gamma-standar}
\begin{array}{l}
b'_n=\theta\big(\log n+ (\nu-1)\log \log n -\log \Gamma(\nu)\big),\\
a'_n=\theta.
\end{array}
\end{equation}
In the proposal of $b''_n$ for a generalized Weibull distribution we
just added one more term  to $b'_n$ coming from the asymptotic expansion
of the solution of the first equation of (\ref{exacta-weibull}). However,
here, that addition does not improves sufficiently the approximation
to the Gumber law, and we need first to enhance the habitual  tail equivalent
 distribution to Gamma law, and later to add an additional term to $b'_n$.
Our proposal is to use:

\bigskip

\noindent $\bullet$ \ If   $\bs{\nu \in (1,2]}$:
$$b''_n=\theta\Big(\log n +(\nu-1)\log \log\big(n/\Gamma(\nu)\big)
-\log \Gamma(\nu)
+\frac{(\nu-1)^2\log \log \big(n/\Gamma(\nu)\big)+\nu-1}{\log \big(n/\Gamma(\nu)\big)}\Big).$$

\bigskip

\noindent $\bullet$ If \ $\bs{\nu> 2}$:
\begin{equation*}
b_n''=\theta\Big(\log n+(\nu-1)\log B_n -\log\Gamma(\nu)
+ \frac{(\nu-1)^2\log B_n-(\nu-1)^2\log(\nu-1)+\nu-1}{B_n}\Big),
\end{equation*}
where
$$B_n=\log n+(\nu-1)\log(\nu-1)-\log\Gamma(\nu).$$

In both cases, we propose
\begin{equation*}
a''_n=\dfrac{b''_n}{ b''_n/\theta-\nu+1}.
\end{equation*}

Note that for $\nu \in (1,2]$,
$$b''_n-b'_n=\theta\Big((\nu-1)\log \frac{\log n-\log \Gamma(\nu)}{\log n}+\frac{(\nu-1)^2\log \log \big(n/\Gamma(\nu)\big)+\nu-1}{\log \big(n/\Gamma(\nu)\big)}\Big).$$
which goes to 0 as $n\to \infty$.
In the other case the difference is similar.

\section{The simplest case}
\label{sec-casfacil}
As we commented, in order to show the problem at hand and  the techniques that we use,
we  first consider
the following particular case of a generalized Weibull distribution:
\begin{equation}
\label{dist-weibull-s}
F(x)=\begin{cases}
1-e\, x\,e^{-x},& \text{if $x\ge 1,$}\\
0, & \text{if $x<1$,}
\end{cases}
\end{equation}
where  $K=e$ and $ C=\tau=\alpha=1$. The standard constants (\ref{aprox-weibull})  are
\begin{equation}
\label{aprox-facil}
\begin{array}{l}
b'_n=\log n+\log \log  n+1,\\
a'_n=1.
\end{array}
\end{equation}
Consider a sample size $n=100$ from the distribution $F$. Solving numerically the first equation of
 (\ref{exacta-weibull}) we get $b_n\approx 7.6384$ (see next Subsection), and from the second equation
 we obtain $a_n\approx 1.1506$. The standard constants are $b'_n\approx 7.1323$ and
 $a'_n=1$. In Figure \ref{densitats}  there is a plot of the density of the Gumbel law
 and the densities of the random
variables $$Y_n=\frac{1}{a_n}(M_n-b_n)
\qquad \text{and} \qquad Y'_n=\frac{1}{a'_n}(M_n-b'_n).$$

\begin{figure}[htb]
\centering
\includegraphics[width=0.6\linewidth]{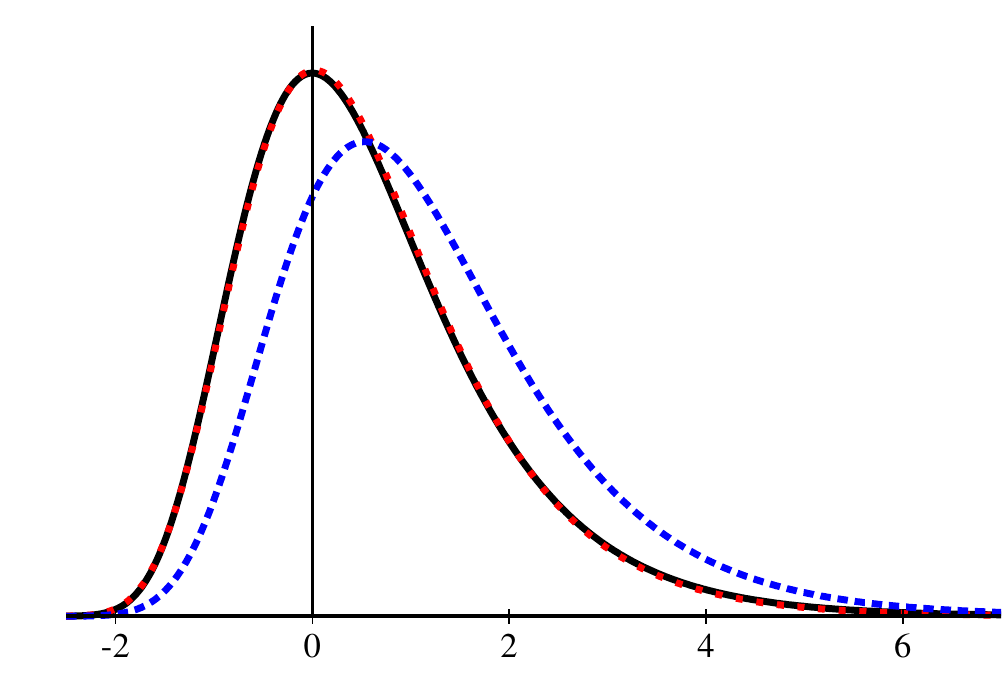}
\caption{ Solid line:
Gumbel density. Red dotted line: Density of $Y_n$. Blue dashed  line:
Density of $Y'_n$.}
\label{densitats}
\end{figure}

In Figure \ref{densitats} we observe that the densities of the Gumbel law and of $Y_n$ are practically
indistinguishable; on the contrary, the density of $Y'_n$ is indeed quite far from
the Gumbel density. As a consequence, this plot illustrates that for $n=100$,
 the distribution of $Y_n$
is very near to the limit, but $Y'_n$ it is not, so
 the approximate norming
constants are very important for such a sample size (and much bigger sample
 sizes, see Table \ref{comparation}).
Then we try  to improve the accuracy of $Y'_n$ choosing
other norming constants.
To get some insight in this question we study  the  first equation of
 (\ref{exacta-weibull}) using the Lambert W function.

\subsection{The  Lambert W function}
\label{sub:lambert}
In the real case  the Lambert function  $W$
is defined implicitly through the real solution of the equation
$$W(x)\, e^{W(x)}=x.$$
Equivalently, $W$ is the inverse of the function  $f(t)=t\,e^t$.  A plot
of the function $f$ (see Figure \ref{lambert-inv})
shows that the Lambert function has two branches (see Figure \ref{lambert-graf}), the   principal one, denoted
by $W_0$, is defined on  $(-1/e,\infty)$, and the secondary, denoted by
$W_{-1}$, is defined on  $(-1/e,0)$ and it satisfies
$$\lim_{x\to 0^-}W_{-1}(x)=-\infty.$$

\begin{figure}[htb]
\centering
\includegraphics[width=0.5\linewidth]{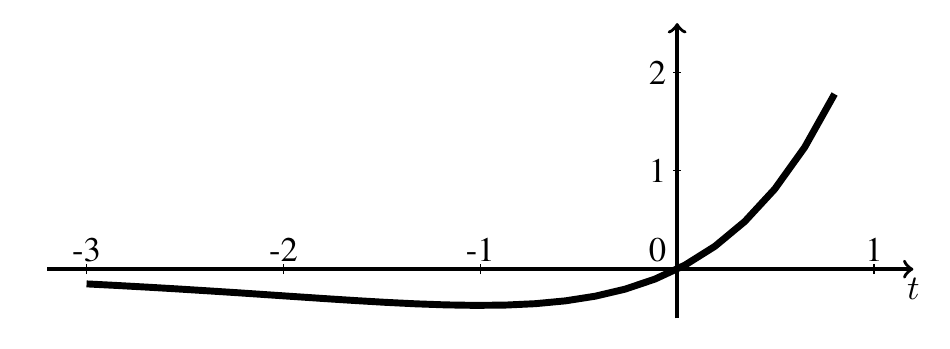}

\caption{Function $f(t)=t\,e^{t}$}
\label{lambert-inv}
\end{figure}

\begin{figure}[htb]
\centering
\includegraphics[width=0.5\linewidth]{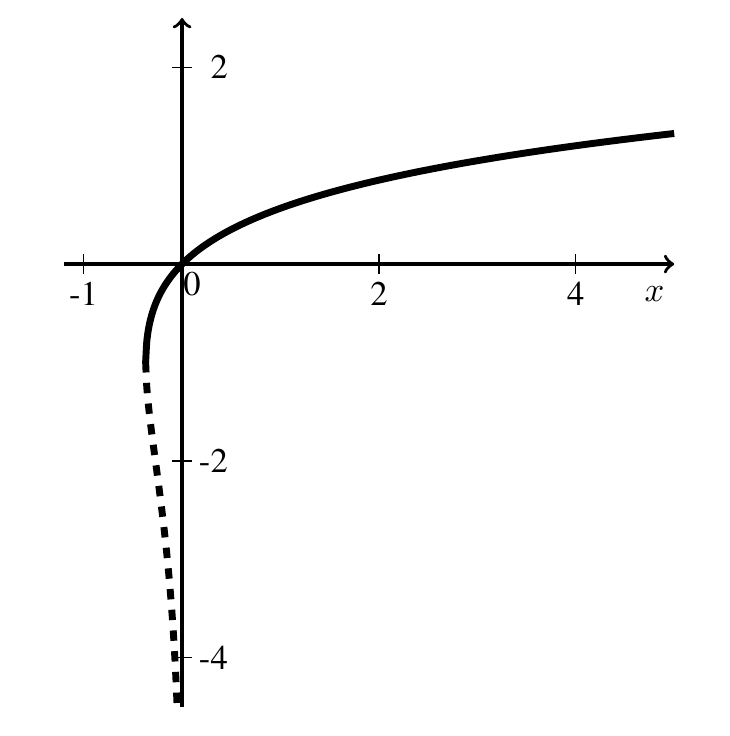}
\caption{Solid line: Principal branch of the real Lambert   W function.
Dashed line: Secondary branch.}
\label{lambert-graf}
\end{figure}

We are interested in the secondary branch (see Subsection \ref{sub:lambert});  its  asymptotic expansion
is
(Corless {\it et al.} \cite[pp. 22 and 23]{CorGonHarJefKnu86},
 see also Comtet \cite{Com70} for the expression of the polynomials and
De Bruijn \cite[pp. 25--27]{DeB81})
\begin{equation}
\label{expansion}
W_{-1}(x)=L_1(-x)-L_2(-x)+\sum_{n=1}^\infty(-1)^{n+1}\frac{P_n(L_2(-x))}{L_1^n(-x)}, \ x\to 0^-,
\end{equation}
where
\begin{equation}
\label{ls}
L_1(x)=\log x \quad\text{and}\quad  L_2(x)=\log \vert \log x\vert,
\end{equation}
and $P_n(x)$ are polynomials related with the signed Stirling numbers
of the first type; the first three polynomials are
$$P_1(x)=x,  \ P_2(x)=\frac{1}{2}\,x^2-x\quad \text{and}\quad P_3(x)=\frac{1}{3}\,x^3-
\frac{3}{2}\,x^2+x.$$
A partial sum approximation to the series on the right hand side of (\ref{expansion})  can be given in the following way
\begin{equation*}
W_{-1}(x)=L_1(-x)-L_2(-x)+\sum_{n=1}^N(-1)^{n+1}\frac{P_n(L_2(-x))}{L_1^n(-x)}+
O\bigg(\Big(\frac{L_2(-x)}{L_1(-x)}\Big)^{N+1}\bigg).
\end{equation*}

\begin{notations}

As usual, we write that $g(x)=O(h(x))$ when $x\to \infty$ if there is
a point  $x_0$ and a constant $C$ such that $\vert g(x)\vert \le Ch(x)$, for all $x>x_0$
(it is assumed $h(x)>0$ for $x>x_0$).
We write  $g(x)=o(h(x))$ if
$\lim_{x\to\infty}\dfrac{g(x)}{h(x)}=0$ (again, here,  $h(x)>0$ for $x>x_0$, for some $x_0$).
Finally, we say that $g(x)$ and $f(x)$ are asymptotically equal and write
$h\sim g$ if $\lim_{x\to\infty}\dfrac{g(x)}{h(x)}=1$. Similar notations are used
when we consider $x\to a$

\end{notations}

\begin{remark}
In this paper, all the computations related with the Lambert W function
are done with the function {\tt lambertW} of the package {\tt emdbook} of
the software {\bf R}.
\end{remark}

\subsection{Computation of the norming constants via  Lambert function}
\label{sub:comp:lambert}
To compute the norming constants,   the first equation of  (\ref{exacta-weibull})
for $F$ given in (\ref{dist-weibull-s}) is
\begin{equation*}
eb_ne^{-b_n}=\frac{1}{n},
\end{equation*}
and hence,
\begin{equation}
\label{branca-in}
b_n=-W\Big(-\frac{1}{e n}\Big).
\end{equation}
By construction
$$\lim _n b_n=\lim_n F^{-1}(1-n^{-1})=\infty,$$
so  in   (\ref{branca-in}) it is needed to consider the secondary branch:
\begin{equation}
\label{exacta-facil-lambert}
b_n=-W_{-1}\Big(-\frac{1}{e n}\Big).
\end{equation}

The asymptotic behaviour of    $b_n$  can be deduced from  (\ref{expansion})  and gives
\begin{align}
\label{desenv-lambert}
b_n=\log  n+1+\log  (\log  n+1) +\frac{\log  (\log  n+1)}{\log  n+1}
+O\Bigg(\bigg(\frac{\log( \log  n+1)}{\log  n+1}\bigg)^2\Bigg).
\end{align}

From the second  equation of  (\ref{exacta-weibull}) we deduce
$$a_n=\frac{b_n}{b_n-1}.$$

\subsection{Computation of the velocity of convergence}

The convergence  (\ref{conv0}) is equivalent that
for every  $x\in \R$,
\begin{equation*}
\label{conv}
\lim_nP\Big\{\frac{1}{a_n}(M_n-b_n)\le x\Big\}=\Lambda(x),
\end{equation*}
or
\begin{equation*}
\label{conv2}
\lim_n F^n (a_nx+b_n)=\Lambda(x).
\end{equation*}
We will prove in Theorem  \ref{velocitat} that
\begin{equation*}
F^n(a_nx+b_n)=\Lambda(x)\Big(1+O\big(1/\log(n)\big)\Big).
\end{equation*}
In this   expression  the constant implicit in $O\big(1/\log n)$ may depend
on $x$.

Furthermore, see again
Theorem  \ref{velocitat},  if $\wt b_n$
satisfies
 $$
 \lim_n\dfrac{b_n-\wt b_n}{a_n}=0,$$
 and
$$\wt a_n=1+O\big(1/\wt b_n\big),$$
then
$$F^n(\wt a_nx+\wt b_n)=\Lambda(x)\Big(1+O\big(1/\log n\big)+ O\big({\wt b_n}/{b_n}-1\big)+O\big(b_n-\wt b_n\big)\Big).$$
In particular, for the standard constants (\ref{aprox-facil}),
\begin{equation}
\label{standard-facil}
b'_n=\log n+\log\log n+1,
\end{equation}
we have
$$b_n-b'_n=O\Big(\frac{\log \log n }{\log n}\Big).$$
However, if we take  one more term of the asymptotic expansion of Lambert W
function in agreement with
(\ref{desenv-lambert}):
\begin{equation}
\label{suggested-facil}
b''_n=\log  n+1+\log  (\log  n+1) +\frac{\log ( \log  n+1)}{\log  n+1}.
\end{equation}
then
$$b_n-b''_n=O\bigg(\Big(\frac{\log \log n }{\log n}\Big)^2\bigg).$$
As it is shown in Table \ref{comparation} the improvement is remarkable, specially for
moderate sample sizes.

\begin{table}[htb]
\centering
\begin{tabular}{rcccccc}
\toprule
{$\bs n$} & $10$  & $10^2$ & $10^3 $ & $10^4 $  & $10^5 $ & $10^6 $ \\
$b_n$ &  4.8897 & 7.6384 &10.2334 &12.7564 & 15.2366& 17.6884
  \\
$b'_n$ & 4.1366 &  7.1323 &  9.8404& 12.4307& 14.9564& 17.4413
\\
$b''_n$ & 4.8590 & 7.6364 & 10.2371 & 12.7613 & 15.2416 & 17.6931
\\
\bottomrule
\end{tabular}
\caption{Comparison of the constants for the distribution
(\ref{dist-weibull-s}):
$b_n$ is computed  numerically, $b'_n$ is  the standard constant (\ref{standard-facil}),
and $b_n''$ the proposed constant (\ref{suggested-facil}).}
\label{comparation}
\end{table}

The  difference $b_n-b'_n$ decreases very slowly, whereas $b_n$ and $b''_n$
are practically indistinguishable; it is worth noting that the difference
between $b'_n$ and $b''_n$ is  a sequence that converges to zero but so slowly
that is important in typical cases. A plot of $b_n,\, b'_n$ and $b''_n$ for values of $n$ from 10 to 1000 is given in Figure \ref{dndprima}.
The lines corresponding to   $b_n$ and $b''_n$ are indistinguishable.
Of course, we can add more terms to $b_n''$, but Table \ref{comparation} suggests
that it is
unnecessary.

\begin{figure}[htb]
\centering
\includegraphics[width=0.5\linewidth]{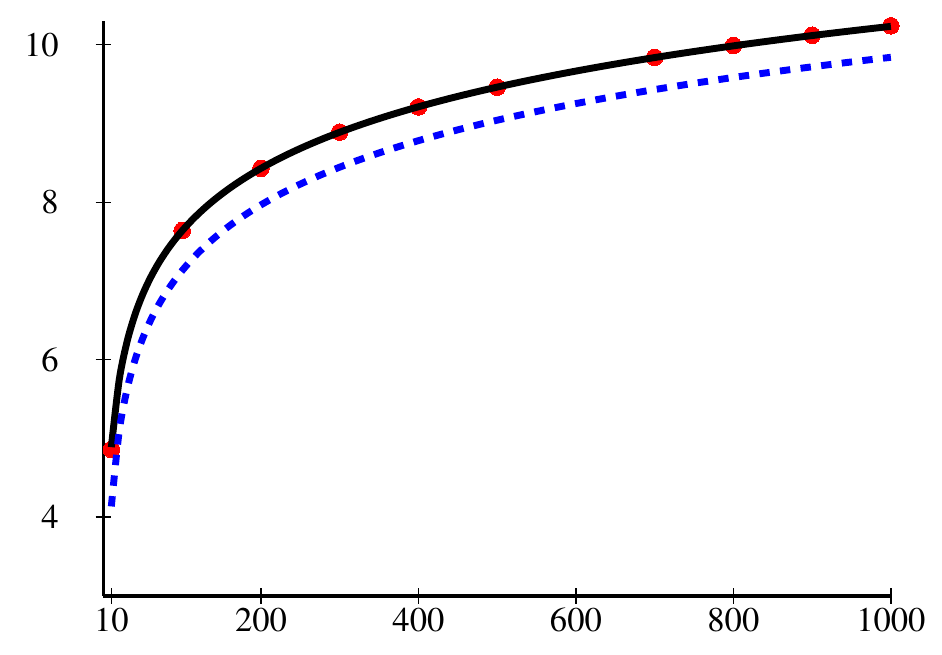}
\caption{Solid line: Numerical value $b_n$. Dashed blue line: standard $b_n'$. Red dots:
 Proposed  $b''_n$.}
\label{dndprima}
\end{figure}

\subsection{The simulation approach}
\label{simulacio}
Random variables with distribution function given by (\ref{dist-weibull-s}) can be
simulated by inversion method because the quantile function corresponding
to $F$   is explicit in terms of the Lambert W function:
$$F^{-1}(u)=-W_{-1}((u-1)/e), \ u\in (0,1).$$
The random variables used to construct the  histograms in Figure \ref{dosHistogrames}
were simulated in this way.

\section{Generalized Weibull distribution}
\label{sec:weibullgeneral}
In this section we deal with a sample  $X_1,\dots,
X_n$ of i.i.d. random variables with generalized  Weibull distribution as presented
in Section \ref{sec:intro}, with distribution function
$$F(x)=1-Kx^\alpha \exp\{-Cx^\tau\},\ x\ge x_0,$$
with
 $\alpha>0 $ and $\tau\ge 1$.
 Let $a_n$ and $b_n$ be given by
 \begin{equation}
\label{exacta-weibull2}
\begin{array}{l}
b_n=F^{-1}(1-n^{-1}),\\
a_n=\dfrac{1}{C\tau b_n^{\tau-1}-{\alpha}/{b_n}}.
\end{array}
\end{equation}

 \begin{theorem}
 \label{velocitat} With the previous notations,
 \begin{enumerate}[\bf 1.]
\item \begin{equation}
\label{velocitat-gen}
F^n(a_nx+b_n)=\Lambda(x)\bigg(1+O\Big(\frac{1}{b_n^\tau}\Big)\bigg)=
\Lambda(x)\Big(1+O\bigg(\frac{1}{\log n}\Big)\bigg),
\end{equation}
where $\Lambda(x)=\exp\{-e^{-x}\}$ is the distribution function of the Gumbel law.

\item
 If $\wt b_n$
satisfies
\begin{equation}
\label{cond.bn}
\lim_n\dfrac{b_n-\wt b_n}{a_n}=0,
\end{equation}
and
\begin{equation}
\label{cond.an}
  \wt a_n=\frac{1}{C\tau \wt b_n^{\tau-1}}
+O\Big(\frac{1}{ \wt b_n^{2\tau-1}}\Big),
\end{equation}
then
$$F^n(\wt a_nx+\wt b_n)=\Lambda(x)\bigg(1+O\big(b_n^\tau-\wt b_n^\tau)+
O\Big(\big(\frac{\wt b_n}{b_n}\big)^\alpha-1\Big)+O\Big(\frac{1}{b_n^\tau}\Big)\bigg).$$

\item
For every sequence of norming constant  $\{\wh a_n,\, n\ge 1\} $
of the form
\begin{equation}
\label{ahat}
 \wh a_n=\frac{1}{C\tau  b_n^{\tau-1}}+\frac{\delta}{ b_n^{2\tau-1}}
+O\Big(\frac{1}{b_n^{2\tau}}\Big),
\end{equation}
we have
$$F^n(\wh  a_nx+b_n)=\Lambda(x)\bigg(1+
\big(2(\alpha-C^2\tau^2 \delta)x-(\tau-1)x^2\big)O\Big(\frac{1}{b_n^{\tau}}\Big)+ O\Big(\frac{1}{b_n^{\tau+1}}\Big)\bigg).$$
As a consequence, when $\tau=1$ and $\delta=\alpha/C^2$,  the sequence $\{a_n,\, n\ge 1\}$ is optimal between
all   sequences $\{\wh a_n,\, n\ge 1\} $ of the above form (\ref{ahat}).

\end{enumerate}
 \end{theorem}

\begin{remarks}
\rule[-1mm]{0cm}{1mm}
\begin{enumerate}[(i)]
\item The constant $a_n$ given in (\ref{exacta-weibull2}) satisfies (\ref{cond.an}) with $b_n$; specifically,
\begin{equation}
\label{desen.an}
a_n=\frac{1}{C\tau b_n^{\tau-1}}\,\Big(1-\frac{\alpha}{C\tau b_n^\tau}\Big)^{-1}=
\frac{1}{C\tau b_n^{\tau-1}}+\frac{\alpha}{C^2\tau^2 b_n^{2\tau-1}}+O\Big(\frac{1}{b_n^{3\tau-1}}\Big).
\end{equation}
Note that it has  the form (\ref{ahat}) with $\delta=\alpha/(C\tau)^2.$

 \item It holds that $\lim_n(b_n^\tau-\wt b_n^\tau)=0$. This is proved in the following
 way: from  (\ref{cond.bn}) and the expression of $a_n$ given in (\ref{exacta-weibull2})
 $$\lim_n\frac{\wt b_n- b_n}{1/b_n^{\tau-1}}=0,$$
 that is,
$$\wt b_n= b_n+ o\Big(\frac{1}{b_n^{\tau-1}}\Big)=b_n\,\bigg(1+o\Big(\frac{1}{b_n^{\tau}}\Big) \bigg).$$
Thus
$$\wt b_n^\tau=b_n^\tau \bigg(1+o\Big(\frac{1}{b_n^{\tau}}\Big) \bigg),$$
and
$$\wt b_n^\tau-b_n^\tau=o(1).$$
\item The case $\tau=1$ is important because a Gamma law is Weibull-like with such a $\tau$
(see Section \ref{chi2}). In this case, in agreement with Remark (i), the sequence
$\{a_n,\, n \ge 1\}$ is optimal.

 \end{enumerate}
\end{remarks}

\bigskip

To prove this theorem we need the following lemma:

\begin{lemma}
\label{lemma}
Let  $A\ne 0$ and consider two sequences
$\{c_n,\,n\ge 1\}$ and $\{d_n,\,n\ge 1\}$ such that
$\lim_n c_n=0$ and $\lim_n d_n=1$
Then, when  $n\to \infty,$
$$ \bigg(1+\frac{A d_n}{n}\big(1+c_n\big)\bigg)^n=e^{A}\Big(1+Ac_n+A(d_n-1)+O(1/n)\Big).$$
\end{lemma}

\noindent{\it Proof of the lemma.}

\smallskip

From  the asymptotic approximation
$$\log(1+x)=x+O(x^2),\ \text{when}\ x\to 0,$$
it follows that
\begin{equation*}
n\log\bigg(1+\frac{Ad_n}{n}\big(1+ c_n\big)\bigg)=
Ad_n(1+c_n)+nO\big(d_n^2(1+c_n)^2/n^2\big)=Ad_n(1+c_n)+O(1/n).
\end{equation*}
Then
\begin{equation*}
 \bigg(1+\frac{Ad_n}{n}\big(1+c_n\big)\bigg)^n=e^{A}e^{A c_n+ A(d_n-1)+O(1/n)}=
 e^{A}\Big(1+Ac_n+A(d_n-1)+O(1/n)\Big).\qquad
 \blacksquare
\end{equation*}

\bigskip

\noindent{\it Proof of Theorem \ref{velocitat}}.

\medskip

\noindent{\bf 1.} In this proof we will use
 that $b_n$ is asymptotically equivalent to $\log^{1/\tau} n$, that
is proved in Subsection \ref{subsec:lam}.
First, note that from (\ref{exacta-weibull2}),
$$a_n=\dfrac{1}{C\tau b_n^{\tau-1}-{\alpha}/{b_n}}=O\Big(\frac{1}{b_n^{\tau-1}}\Big),$$
and
$$\frac{a_n}{b_n}=O\Big(\frac{1}{b_n^{\tau}}\Big).$$
Then
\begin{align*}
K(a_nx+&b_n)^\alpha  \exp\{-C(a_nx+b_n)^\tau\}= Kb_n^\alpha
\exp\{-Cb_n^\tau\} \Big(\frac{a_n}{b_n}x+1\Big)^\alpha
\exp\{-C(a_nx+b_n)^\tau+Cb_n^\tau\} \notag \\
&\buil{=}{}{(*)} \frac{e^{-x}}{n} \Big(1+O\big(\frac{a_n}{b_n}\big)\Big)
\exp\Big\{-Cb_n^\tau\Big(\big(\frac{a_n}{b_n} x+1\big)^\tau-1\Big)+x\Big\}\notag \\
&=\frac{e^{-x}}{n} \Big(1+O\Big(\frac{1}{b_n^{\tau}}\Big)\Big)
\exp\bigg\{-Cb_n^\tau\bigg(1+\tau \frac{a_n}{b_n}x+O\Big(\big(\frac{a_n}{b_n}\big)^2\Big)-1\bigg)+x\bigg\} \notag\\
&=\frac{e^{-x}}{n} \Big(1+O\Big(\frac{1}{b_n^{\tau}}\Big)\Big)
\exp\bigg\{\Big(-C\tau  b_n^{\tau-1}a_n x+O\Big(\frac{1}{b_n^\tau}\Big)+x\bigg\} \notag\\
&=\frac{e^{-x}}{n} \Big(1+O\Big(\frac{1}{b_n^{\tau}}\Big)\Big)
\exp\bigg\{\Big(-\frac{C\tau  b_n^{\tau}}{C\tau b_n^{\tau}-\alpha}\, x+O\Big(\frac{1}{b_n^\tau}\Big)+x\bigg\} \notag\\
&=\frac{e^{-x}}{n} \Big(1+O\Big(\frac{1}{b_n^{\tau}}\Big)\Big)
\exp\Big\{O\Big(\frac{1}{b_n^{\tau}}\Big)\Big\}
=\frac{e^{-x}}{n}\bigg(1+O\Big(\frac{1}{b_n^{\tau}}\Big)\bigg),
\end{align*}
where the equality $(^*)$ follows from the definition of $b_n$ given
in (\ref{exacta-weibull2}).
Since $\lim_n 1/b_n=0$,   we
can apply Lemma \ref{lemma} (with $d_n=1$):
$$F^n(a_nx+b)=\Big(1-K(a_nx+b_n)^\alpha  \exp\{-C(a_nx+b_n)^\tau\}\Big)^n
=\exp\{-e^{-x}\}\Big(1+O\big(\frac{1}{b_n^{\tau}}\big)+O\big(\frac{1}{n}\big)\Big),$$
and from the estimation $b_n\sim
\log^{1/\tau} n$, the term $O(1/n)$ can be eliminated.

\bigskip

\noindent{\bf 2.}

Proceeding as before,
\begin{align}
\label{dif-aprox}
K(\wt a_nx+&\wt b_n)^\alpha  \exp\big\{-C(\wt a_nx+\wt b_n)^\tau\big\}\notag\\
&=Kb_n^\alpha
\exp\{-Cb_n^\tau\} e^{-x}\Big(\frac{\wt a_nx+\wt b_n}{b_n}\Big)^\alpha
\exp\big\{-C(\wt a_nx+\wt b_n)^\tau+C b_n^\tau+x\big\}\notag\\
&= \frac{e^{-x}}{n}\Big(\frac{\wt a_nx+\wt b_n}{b_n}\Big)^\alpha
\exp\big\{-C(\wt a_nx+\wt b_n)^\tau+C b_n^\tau+x\big\}
\end{align}

The term in the exponential is 
\begin{align*}
-C(\wt a_nx+& \wt b_n)^\tau+C b_n^\tau+x
=-C \wt b_n^\tau \Big(1+\frac{\wt a_n}{\wt b_n}\, x\Big)^{\tau}+C b_n^\tau+x\\
&=-C \wt b_n^\tau \bigg(1+\tau\,\frac{\wt a_n}{\wt b_n}\, x
+O\Big(\big(\frac{\wt a_n}{\wt b_n}\big)^2\Big) \bigg)+C b_n^\tau+x\\
&=-C \big(\wt b_n^\tau-b_n^\tau\big)- C\tau \wt a_n\wt b_n^{\tau-1}
x + x +O\big(\frac{1}{\wt b_n^{\tau}}\big)\\
&=-C \big(\wt b_n^\tau-b_n^\tau\big)+O\big(\frac{1}{\wt b_n^{\tau}}\big),
\end{align*}
where in the last equality we used (\ref{cond.an}).

Then
\begin{align*}
(\ref{dif-aprox})&= \frac{e^{-x}}{n}\Big(\frac{\wt b_n}{b_n}\Big)^\alpha
\,\Big(1+\frac{\wt a_n}{\wt b_n}\, x\Big)^\alpha\,
\Big(1-C \big(\wt b_n^\tau-b_n^\tau\big)+O\big(\frac{1}{\wt b_n^{\tau}}\big))\\
&=\frac{e^{-x}}{n}\Big(\frac{\wt b_n}{b_n}\Big)^\alpha
\Big(1-C \big(\wt b_n^\tau-b_n^\tau\big)+O\big(\frac{1}{\wt b_n^{\tau}}\big)\Big)\\
\end{align*}
and we apply again Lemma  \ref{lemma}.

\bigskip

\noindent{\bf 3.} Developing  as in {\bf 1}, but now  taking one more term in the
expansions
\begin{align*}
K(\wh  a_n & x+b_n)^\alpha  \exp\{-C(\wh a_nx+b_n)^\tau\}=
\frac{e^{-x}}{n}\Big(\frac{\wh a_n}{b_n}x+1\Big)^\alpha
\exp\Big\{- C b_n^\tau\Big(\frac{\wh a_n}{b_n}\, x+1\Big)^\tau+Cb_n^\tau+x\big\} \notag \\
&= \frac{e^{-x}}{n} \bigg(1+\alpha\, \frac{\wh a_n}{b_n}\, x+
O\Big(\big(\frac{\wh a_n}{b_n}\big)^2\Big)\bigg)\\
&\qquad \times\exp\Big\{-Cb_n^\tau\Big(1+\tau\,\frac{\wh a_n}{b_n}\, x+
 \frac{\tau(\tau-1)}{2}\,x^2
\big(\frac{\wh a_n}{b_n}\big)^2+O\Big(\big(\frac{\wh a_n}{b_n}\big)^3\Big)
-1\Big)+x\Big\}\notag \\
&=\frac{e^{-x}}{n} \bigg(1+\frac{\alpha}{C\tau b_n^\tau}\, x+
O\Big(\frac{1}{b_n^{2\tau}}\Big)\bigg)\,\exp\Big\{-\frac{C\tau\delta}{b_n^\tau}
\, x- \frac{\tau-1}{2 C\tau b_n^\tau}\, x^2+
O\Big(\frac{1}{b_n^{\tau+1}}\Big)\Big\}\\
&=\frac{e^{-x}}{n} \bigg(1+\frac{2\big(\alpha-C^2\tau^2\delta\big) x-(\tau-1)x^2}{2C\tau}
\,\frac{1}{b_n^\tau}+
O\Big(\frac{1}{b_n^{\tau+1}}\Big)\bigg).
\end{align*}
By (\ref{desen.an}), when $\wh a_n=a_n$, then $\delta=\alpha/(C\tau)^{2},$
and if, moreover, $\tau=1$, all terms of order $b_n^{-\tau}$ cancel. \qquad $\blacksquare$

\subsection{Computation of the constants using Lambert W function}
\label{subsec:lam}
The first equation of (\ref{exacta-weibull2}) in this case is
\begin{equation}
\label{eq-weibull-general}
Kb_n^\alpha \exp\{-Cb_n^\tau\}=\frac{1}{n}.
\end{equation}
In terms of the  Lambert W function  the  solution is
$$b_n=\Bigg(-\frac{\alpha}{C\tau}\,W_{-1}\Big(-\frac{C\tau}{\alpha (Kn)^{\tau/\alpha}}\Big)\Bigg)^{1/\tau}.$$
From the asymptotic expansion of the Lambert function (\ref{expansion}), we have that
\begin{equation}
\label{general-dn-exp}
b_n=\Big(\frac{\alpha}{C\tau}\Big)^{1/\tau}\bigg(-M_1+M_2-\frac{M_2}{M_1}+\cdots\bigg)^{1/\tau}
\end{equation}
where,
\begin{align*}
M_1&=L_1\Big(\frac{C\tau}{\alpha (Kn)^{\tau/\alpha}}\Big)=
\log\Big(\frac{C\tau}{\alpha (Kn)^{\tau/\alpha}}\Big),\\
M_2&=L_2\Big(\frac{C\tau}{\alpha (Kn)^{\tau/\alpha}}\Big)=\log (-M_1).
\end{align*}

\subsection{Computation of the constants using  Comtet expansion}
\label{subsec:comtet}
The solution of equation (\ref{eq-weibull-general})
can be expressed in an alternative way:
Fix $\gamma\ne 0$ and denote by $U_\gamma(x)$  the (unique) solution
$t$ of the equation
$$t^\gamma e^t=x$$
such that $t\to\infty$ when $x\to \infty$.
Equation
(\ref{eq-weibull-general}) is equivalent to
$$\frac{1}{K}\, b_n^{-\alpha}e^{Cb_n^\tau}=n,$$
or
$$\big(Cb_n^\tau\big)^{-\alpha/\tau} e^{Cb_n^\tau}=C^{-\alpha/\tau}Kn.$$
Hence,
\begin{equation}
\label{wei-comtet}
b_n=\bigg(\frac{1}{C}\,U_{-\alpha/\tau}\Big(\frac{Kn}{C^{\alpha/\tau}}\Big)\bigg)^{1/\tau}.
\end{equation}

Comtet \cite{Com70} extended De Bruijn
 expansion of the principal branch of  the Lambert W function
to $U_\gamma$ obtaining
\begin{equation}
\label{expansion-comtet}
U_\gamma(x)=L_1(x)-\gamma L_2(x)+\sum_{n=1}^\infty(-\gamma)^{n+1}\frac{P_n(L_2(x))}{L_1^n(x)},
\ x\to\infty,
\end{equation}
where, $L_1(x), \, L_2(x)$ and $P_n(x)$ are as in   (\ref{expansion}) and (\ref{ls}).
So applying that expansion to (\ref{wei-comtet}) we get a new asymptotic expansion
for $b_n$. The finite expansion obtained by truncation of  (\ref{expansion-comtet}) and
the one corresponding to (\ref{general-dn-exp}) are not equal;
the difference is originated from the fact that some constants
appearing early in
(\ref{general-dn-exp}) are delayed in (\ref{wei-comtet}): the whole sum of both
series is the same, but the truncated series are slightly different.
It turns out that
when $\alpha>\tau$ then it is better the truncation from (\ref{general-dn-exp}),
when $\alpha<\tau$ it is better the approximation given by Comtet expansion
(\ref{expansion-comtet}), and when $\alpha=\tau$ both expansion coincide. See Appendix A1
for a proof.

\subsection{Return to the norming constants}
Returning to the norming constants of the general Weibull maxima, when
 $\alpha> \tau$,
the better finite  asymptotic expansion is given by the truncation of   (\ref{general-dn-exp}),
and we get the constants given in (\ref{bsecgeneral}).
When $\alpha\le  \tau$, we use the asymptotic expansion (\ref{expansion-comtet})
and we deduce the constant given in (\ref{bsecalphapetit}).

\section{Maxima of Gamma random variables}
\label{chi2}
In this section we study the case  when the sample comes from a Gamma distribution.
That case appears (for $\chi^2(m)$ laws) in  some practical problems of  signal analysis
 and, as we commented in the Introduction, it  is in the origin of this work.
Since a $\chi^2(m)$  law is a Gamma$(\nu,\theta)$ with $\nu=m/2$ and $\theta=2$, we study
directly the latter case.
A Gamma law is Weibull--like with $\alpha=\nu-1$ and $\tau=1$ (see below); we will also
 assume that $\nu>1$ (in agreement with our general restriction $\alpha>0$).

The distribution function $G(x)$ of a  Gamma($\nu,\theta$) law  can be written in terms of the incomplete
Gamma function as
\begin{equation}
\label{gamma}
G(x)=1-\frac{\Gamma(\nu,x/\theta)}{\Gamma(\nu)}, \ x>0,
\end{equation}
where $\Gamma(a)$  is the  Gamma function,
$$\Gamma(a)=\int_0^\infty t^{a-1}e^{-t}\, dt,\ a>0,$$
and $\Gamma(a,y)$ is the  upper  incomplete Gamma function
$$\Gamma(a,y)=\int_y^\infty t^{a-1}e^{-t}\, dt,\ a,y>0$$
Then,  from the asymptotics for the incomplete
Gamma function deduced from Olver {\it et al.} \cite[formula 8.11.2]{NIST10}

$$\lim_{y\to\infty}\frac{\Gamma(a,y)}{y^{a-1}e^{-y}}=1,$$ it follows
\begin{equation}
\label{weibull-approx}
\lim_{x\to \infty}\frac{1-G(x)}{Kx^\alpha \exp\{-C x^\tau\}}=1,
\end{equation}
with $$K=\frac{1}{\theta^{\nu-1}\Gamma(\nu)},\quad \alpha=\nu-1,\quad C=\frac{1}{\theta}, \quad\text{and}\quad
\tau=1.$$
So we can consider the tail equivalent function
$$F_1(x)=1-\frac{1}{\theta^{\nu-1}\Gamma(\nu)}\,x^{\nu-1}\exp\{-x/\theta\}, \ x\ge x_0.$$
The  auxiliary function (\ref{aux-canon}) corresponding   to $F_1$ is
\begin{equation}
\label{canonic-wibull}
A_1(x)=\dfrac{x}{x/\theta-\nu+1}.
\end{equation}
This means that the Gamma($\nu,\theta$) law belongs to the Weibull--like distributions.
 In agreement with our comments on Section \ref{sec:intro},
the standard constants (see (\ref{aprox-weibull})) are
\begin{equation}
\label{aprox-chi}
\begin{array}{l}
b'_n=\theta\big(\log n+ (\nu-1)\log (\log n) -\log \Gamma(\nu)\big),\\
a'_n=\theta.
\end{array}
\end{equation}
Numerical computations show that these constants produce very inaccurate
results (see Subsection \ref{sec:constant-chi}). However, in  this case, the simple addition of more terms using
Comtet or Lambert expansion as in previous section does not improve
  the
results
and we need to consider a right tail equivalent distribution function more
  accurate  than $F_1$.
We will use that the incomplete Gamma function $\Gamma(a,y)$ ($a>0$) admits
the following asymptotic  expansion
for $y\to \infty$ (see Olver {\it et al.} \cite[formula 8.11.2]{NIST10}):
\begin{equation*}
\Gamma(a,y) \sim y^{a-1}e^{-y}\bigg(1+\frac{a-1}{y}+\frac{(a-1)(a-2)}{y^2}+\cdots\bigg).
\end{equation*}
Observe that when $a$ is an integer, the series in the right hand side terminates,
and the expression is not only asymptotic but exact for all $y>0$; this is what
happens with the distribution function of a $\chi^2(m)$ random variable with
$m$ even. So, we will consider a distribution function of the form
\begin{equation}
\label{f1}
F_2(x)=1- \frac{1}{\theta^{\nu-1}\Gamma(\nu)}\, x^{\nu-1} \exp\{-x/\theta\}\bigg(1+\frac{\theta(\nu-1)}{x}\bigg),\, x\ge x_0.
\end{equation}

\bigskip

 Now we also need an extension of Lambert and Comtet asymptotic expansions
to this new context.
\subsection{Extension of  Comtet expansion}
\label{comtet-extension:subsec}
Robin \cite{Rob88} and Salvi \cite{Sal92} extended Comtet \cite{Com70} results
in order to deduce an asymptotic expansion of the solution  of the equation
\begin{equation}
\label{eq:salvi}
t^\gamma e^{t} D\Big(\frac{1}{t}\Big)=x,
\end{equation}
such that $t\to \infty$ when $x\to \infty$, where $\gamma\ne 0$ and
$$D(t)=\sum_{n=0}^\infty d_n t^n, \ \text{with} \ d_0\ne 0,$$
is a power series convergent in a neighborhood of the origin. Denote by
 $U_{\gamma,D}(x)$ that solution.
Robin \cite{Rob88} and Salvi \cite{Sal92} prove that for every $N$,
for $x\to\infty,$
\begin{equation}
\label{ugamma}
U_{\gamma, D}(x)=L_1(x)+\sum_{n=0}^N\frac{Q_n(L_2(x))}{L_1^n(x)}
+o\bigg(\frac{1}{L_1^N(x)}\bigg),
\end{equation}
where $L_1$ and $L_2$ are the same as in (\ref{expansion}), and
 $Q_n(x)$ are polynomials
now depending on $\gamma$ and $D$:
$$Q_n=Q_n(\gamma,d_0,\dots,d_n),$$
with degree $n$ for $n\ge 1$, and $Q_0$ has degree 1.
The first two polynomials (fortunately, the only ones that we  need), see Appendix A2,  are
\begin{equation}
\label{q2}
Q_0(x)=-\gamma x-\log d_0  \quad\text{and}\quad  Q_1(x)=\gamma^2 x+\gamma \log d_0-\frac{d_1}{d_0}.
\end{equation}
 When $D(x)=1$, then $Q_0(x)=-\gamma x$, and for $n\ge 1,$
 $Q_n(x)=(-\gamma)^{n+1} P_n(x)$, where
$P_n$ are the Comtet polynomials  in Subsection \ref{sub:lambert}.

\subsection{Extension of  Lambert  expansion}
\label{lambert-extension:subsec}
Consider the equation (\ref{eq:salvi}) for $\gamma=1$.  The inverse of the function
$$f(t)=t e^t D\Big(\frac{1}{t}\Big)$$
($t$ out of a neighborhood of the origen) has a secondary branch,
denoted by $W_{-1,D}(x)$, that goes to $-\infty$ when $x\to 0^-$.
The asymptotic expansion of this branch
(see Appendix A3) is
\begin{align}
\label{lambert-A-extension}
W_{-1,D}(x)=
L_1(-x)+\sum_{n=0}^\infty (-1)^{n+1}\frac{R_n(L_2(-x))}{L_1^n(-x)}.
\end{align}
where the first two polynomials are
\begin{equation}
\label{r2}
R_0(x)=x+\log d_0 \quad \text{and}\quad R_1(x)=x+\log d_0-\frac{d_1}{d_0}.
\end{equation}
The relationship between the polynomials $Q_n$ of section \ref{comtet-extension:subsec}
and $R_n$ is studied in Appendix A4.

\subsection{New norming constants for the  maxima of Gamma random variables}
\label{sec:constant-chi}
We apply the principle of simplification by tail equivalence given in Property
\ref{propietat1},
and from  (\ref{exacta-weibull}) for $F_2$ in (\ref{f1}), $b_n$ verifies
\begin{equation}
\label{dn:gamma}
\Big(\frac{b_n}{\theta}\Big)^{\nu-1} \exp\{-b_n/\theta\}\bigg(1+\frac{\theta(\nu-1)}{b_n}\bigg)=\frac{\Gamma(\nu)}{n}.
\end{equation}
As a consequence of the previous two subsections, we have two ways to express the solution
of this equation.
For  the first one, write $y=b_n/\theta$:  we need to solve
 \begin{equation}
 \label{lambert.eq}
y^{\nu-1}e^{-y}\bigg(1+\frac{\nu-1}{y}\bigg)
=\frac{\Gamma(\nu)}{n},
\end{equation}
 or equivalently,
 \begin{equation*}
y^{1-\nu}e^{y}\bigg(1+\frac{\nu-1}{y}\bigg)^{-1}=\frac{n}{\Gamma(\nu)},
 \end{equation*}
 Hence, with the notation of Subsection \ref{comtet-extension:subsec},
\begin{equation}
\label{dn:comtet-extension}
b_n=\theta \, U_{1-\nu,D}\Big(\frac{n}{\Gamma(\nu)}\Big),
\end{equation}
where,  $D(t)$ is the series
\begin{equation}
\label{D-gamma}
D(t)=\big(1+(\nu-1) t)^{-1}=1-(\nu-1)t+O(t^2).
\end{equation}
In a similar way, we can transform equation (\ref{lambert.eq}):
 \begin{equation*}
y e^{-y/(\nu-1)}\bigg(1+\frac{\nu-1}{y}\bigg)^{1/(\nu-1)}=\Big(\frac{\Gamma(\nu)}{n}\Big)^{1/(\nu-1)},
 \end{equation*}
or
 \begin{equation*}
\Big(-\frac{y}{\nu-1}\Big)\, e^{-y/(\nu-1)}
\bigg(1-\frac{1}{-y/(\nu-1)}\bigg)^{1/(\nu-1)}=
-\frac{1}{\nu-1}\,\Big(\frac{\Gamma(\nu)}{n}\Big)^{1/(\nu-1)}.
 \end{equation*}
Thus, from  Subsection \ref{lambert-extension:subsec},
\begin{equation}
\label{dn:lambert-extension}
b_n=-(\nu-1)\theta\, W_{-1,E}\bigg(-\frac{1}{\nu-1}\,\Big(\frac{\Gamma(\nu)}{n}\Big)^{1/(\nu-1)}\bigg),
\end{equation}
where $E(t)$ is
\begin{equation}
\label{D-gamma2}
E(t)=(1- t)^{1/(\nu-1)}=1-\frac{1}{\nu-1}\, t+O(t^2).
\end{equation}
(We need only the first two terms of this series).
We prove in Appendix A4 that
 the finite  expansions deduced from
(\ref{dn:comtet-extension}) and (\ref{dn:lambert-extension}) when $\nu=2$ are equal,
when $\nu>2$ then (\ref{dn:lambert-extension}) gives a more accurate value,
and when $1<\nu<2$ it is better to use the expansion deduced from
(\ref{dn:comtet-extension}). In both  cases, the results on Section \ref{sec-casfacil}
suggest to add just one term of the asymptotic expansion to the standard constants,
and hence it suffices to consider the polynomials $Q_0$ and $Q_1$, that depend only
on the terms $d_0$ and $d_1$ of the series $D(t)$ or $B(t)$.
The formulas are the following:

\bigskip

{\bf For}  $\bs{\nu \in (1,2]}$: From (\ref{dn:comtet-extension}), (\ref{ugamma})
 and (\ref{q2}), we propose
$$b''_n=\theta\Big(\log\big(n/\Gamma(\nu)\big)+(\nu-1)\log \log\big(n/\Gamma(\nu)\big)
+\frac{(\nu-1)^2\log \log \big(n/\Gamma(\nu)\big)+\nu-1}{\log \big(n/\Gamma(\nu)\big)}\Big).$$
and
\begin{equation}
\label{csec}
a''_n=\dfrac{b''_n}{ b''_n/\theta-\nu+1}.
\end{equation}

\bigskip

{\bf For}  $\bs{\nu\ge 2}$: From (\ref{dn:lambert-extension}), (\ref{lambert-A-extension})
and (\ref{r2}),
\begin{equation}
\label{dsec}
b_n''=\theta\Big(\log n+(\nu-1)\log B_n -\log\Gamma(\nu)
+ \frac{(\nu-1)^2\log B_n-(\nu-1)^2\log(\nu-1)+\nu-1}{B_n}\Big),
\end{equation}
where
$$B_n=\log n+(\nu-1)\log(\nu-1)-\log\Gamma(\nu),$$
and $a_n''$ the same as (\ref{csec}).

\bigskip

In Table \ref{comparison_chi} there are  some numerical results
 for a $\chi^2(10)$ distribution; the numeric  value of $b_n$ is computed
using the quantile function of a $\chi^2$ distribution implemented
in $\bf R$ program.  Similar results
are obtained for the case $\nu\in (1,2)$, for example, for a $\chi^2(3)$ distribution;
however, in this case the discrepancy between the approximation using the standard
constants and the Gumbel distribution is not as grave as in the case $\nu>2$.

\begin{table}[htb]
\centering
\begin{tabular}{rcccccc}
\toprule
{$\bs n$} & $10$  & $10^2$ & $10^3 $ & $10^4 $  & $10^5 $ & $10^6 $ \\
$b_n$ & 15.9872 & 23.2093 & 29.5883 & 35.5640 & 41.2962 & 46.8630
\\
$b'_n$ &   4.9213 & 15.0717 & 22.9606 & 29.8272 & 36.2175 & 42.2812
\\
$b''_n$ &  13.3518 & 22.0874 & 29.0421 & 35.2855 & 41.1581 & 46.8045\\
\bottomrule
\end{tabular}
\caption{Comparison of the constants $b_n$ for the $\chi^2(10)$ distribution:
$b_n$ is the numeric value, $b_n'$ is the standard value, $b''_n$ is the value
given in (\ref{dsec}).}
\label{comparison_chi}
\end{table}
For the norming constant  $a_n$ there are also important differences, see
Table \ref{comparison_c}.

\begin{table}[htb]
\centering
\begin{tabular}{rcccccc}
\toprule
{$\bs n$} & $10$  & $10^2$ & $10^3 $ & $10^4 $  & $10^5 $ & $10^6 $ \\

$a_n$ &  4.0032 & 3.0520 & 2.7411 & 2.5805  & 2.4805 & 2.4117 \\
$a'_n$ &   2    & 2      &   2    &     2  &  2      & 2 \\
$a''_n$ & 4.9896 & 3.1358 & 2.7604 & 2.5864 & 2.4825 & 2.4123 \\
\bottomrule
\end{tabular}
\caption{Comparison of the constants $a_n$ for the $\chi^2(10)$ distribution:
$a_n$ is the numeric value, $a_n'$ is the standard value, $a''_n$ is the value
given in (\ref{csec}).}
\label{comparison_c}
\end{table}

In Figure \ref{densitatschi100} there is a plot of the density funcions of the random
variables
$$Y_n=\frac{1}{a_n}\,(M_n-b_n), \quad   Y'_n=\frac{1}{a'_n}\,(M_n-b'_n),\
 \text{and} \quad Y''_n=\frac{1}{a''_n}\,(M_n-b''_n)$$
from a sample of size  $n=100$ of $\chi^2(10)$ random variables,
where $b_n$ and $a_n$ are the numeric solutions of equations (\ref{exacta-weibull}),
$b'_n$ and $a'_n$ are the standard solution given in (\ref{aprox-chi}),
$b''_n$ and $a''_n$ are the constants 
(\ref{dsec}) and (\ref{csec}).

\begin{figure}[htb]
\centering
\includegraphics[width=0.5\linewidth]{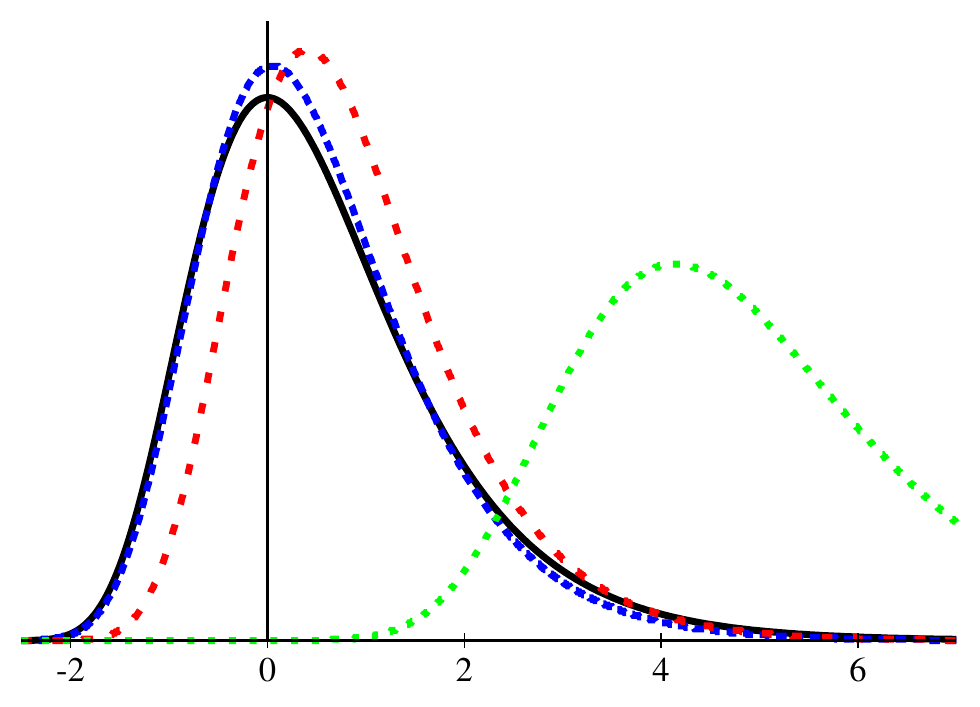}

\caption{ Maximum of 100 $\chi^2(10)$ random variables. Solid line:
Gumbel density. Dashed Blue line: Density of $Y_n$. Loosely dashed red line:
 density of $Y''_n$.
 Dotted green line: Density
of  $Y'_n$.} 
\label{densitatschi100}
\end{figure}

\section*{Appendix}
\subsection*{Appendix A1. Lambert versus Comtet asymptotic expansions}

Fix $\beta > 0$ and consider the solution $t$ of the equation
\begin{equation}
\label{eq:appen}
t^\beta e^{-t}=x,
\end{equation}
such that $t\to \infty$ when $x\to 0^+$.
With the  notations of subsections \ref{subsec:lam} and \ref{subsec:comtet} , the solution can be  written in two ways:
\begin{equation}
\label{dues_expressions}
t=-\beta \,W_{-1}\big(-x^{1/\beta}/\beta\big)=U_{-\beta}(1/x).
\end{equation}
From (\ref{expansion}), the central term of (\ref{dues_expressions}) is
\begin{align}
\label{t-lambert}
t&=-\beta L_1\big(x^{1/\beta}/\beta\big)+
\beta L_2\big(x^{1/\beta}/\beta\big)
-\beta\sum_{n=1}^\infty(-1)^{n+1}\frac{P_n\big(L_2(x^{1/\beta}/\beta\big))}
{L_1^n\big(x^{1/\beta}/\beta\big)}\notag\\
&=-L_1\big(x/\beta^\beta\big)+
\beta L_2\big(x^{1/\beta}/\beta\big)
-\sum_{n=1}^\infty(-\beta)^{n+1}\frac{P_n\big(L_2(x^{1/\beta}/\beta)\big)}
{L_1^n\big(x/\beta^\beta\big)},
\end{align}
where $L_1$ and $L_2$ are defined in (\ref{ls}). From (\ref{expansion-comtet}), the term of the right side of (\ref{dues_expressions})
is
\begin{align}
\label{t-comtet}
t=L_1(1/x)+\beta L_2(1/x)+\sum_{n=1}^\infty \beta^{n+1}\frac{P_n(L_2(1/x))}{L_1^n(1/x)}
\notag\\
=-L_1(x)+\beta L_2(x)-\sum_{n=1}^\infty (-\beta)^{n+1}\frac{P_n(L_2(x))}{L_1^n(x)}.
\end{align}
Comparing (\ref{t-lambert}) and (\ref{t-comtet}) we realize that both asymptotic expansion are
the same function applied to different points; specifically, define
\begin{equation}
\label{function-h}
h(y,z)=-y+\beta z-\sum_{n=1}^\infty (-\beta)^{n+1}\frac{P_n(z)}{y^n}.
\end{equation}
Then (\ref{t-lambert})  is the function $h(y,z)$ applied to
$$y=L_1\big(x/\beta^\beta\big)\quad\text{and}\quad
z=L_2\big(x^{1/\beta}/\beta\big),$$
whereas (\ref{t-comtet})  is the function $h(y,z)$  applied to
$y=L_1(x)$ and $z=L_2(x)$.
The argument of De Bruijn \cite[p. 25--27]{DeB81} in this case shows that
there exists constants $a$ and $b$ such that if $y>a$ and $0<z/y<b$, then the series
in the right hand side of (\ref{function-h}) is absolutely convergent and
$$\big\vert h(y,z)+y-\beta z+\sum_{n=1}^N (-\beta)^{n+1}\frac{P_n(z)}{y^n}\big\vert
\le C\bigg( \frac{z}{y}\bigg)^{N+1}.$$
So we deduce
\begin{equation}
\label{fita-lambert}
\bigg\vert t +L_1\big(x/\beta^\beta\big)-
\beta L_2\big(x^{1/\beta}/\beta\big)
+\sum_{n=1}^N(-\beta)^{n+1}\frac{P_n\big(L_2(x^{1/\beta}/\beta)\big)}
{L_1^n\big(x/\beta^\beta\big)}\bigg\vert \le C \bigg(\frac{L_2\big(x^{1/\beta}/\beta\big)}{L_1\big(x/\beta^\beta\big)}\bigg)^{N+1},
\end{equation}
and
\begin{equation}
\label{fita-comptet}
\bigg\vert t + L_1(x)-\beta L_2(x)+\sum_{n=1}^N (-\beta)^{n+1}\frac{P_n(L_2(x))}{L_1^n(x)}
\bigg \vert \le  C \bigg(\frac{L_2(x)}{L_1(x)}\bigg)^{N+1}.
\end{equation}
Now,
\begin{align}
\label{quocient}
\frac{\dfrac{L_2\big(x^{1/\beta}/\beta\big)}{L_1\big(x/\beta^\beta\big)}}
{\dfrac{L_2(x)}{L_1(x)}}& =
\frac{L_2\big(x^{1/\beta}/\beta\big)}{L_2(x) }
\ \frac{L_1(x)}{L_1\big(x/\beta^\beta\big)}
=\bigg(1-\frac{\log \beta}{L_2(x)}-\frac{\beta\, \log \beta}{L_1(x)}+\cdots\bigg)
\bigg(1+\frac{\beta \log \beta}{L_1(x)}+\cdots\bigg)\notag\\
&=\bigg(1-\frac{\log \beta}{L_2(x)}+\cdots\bigg)
\end{align}
and, for $x>0$ small enough, (remember $L_2(x)=\log(-\log x)<0$, when $x\to 0^+$),
$$(\ref{quocient})\ \text{is}\ \quad\begin{cases}
>1,& \text{if $0<\beta <1$},\\
=1,& \text{if $\beta =1$},\\
<1, & \text{if $\beta>1$}.
\end{cases}
$$

This indicates that

\begin{itemize}
\item If $0<\beta <1$ then the finite  expansion deduced from (\ref{t-comtet})
seems to produce a  more accurate approximation.
\item If $\beta>1$, the finite expansion from (\ref{t-lambert}) seems better.
\item If $\beta =1$, both expansions are equal.
\end{itemize}

Intuitively, (\ref{t-lambert}) and (\ref{t-comtet}) are  asymptotic expansions
when $x\to 0^+$, and, for example,  when $\beta>1$, the dominant part of both, $\log (\dots)$ is applied in
(\ref{t-lambert}) to a smaller number, $x/\beta^\beta$.

In table \ref{comparison-t} there is a  numerical study to illustrate this point.
 We denote by $t$
the numeric solution of equation (\ref{eq:appen}), by $t_W$ the approximation
deduced from Lambert expansion (\ref{t-lambert}),
\begin{equation}
\label{t-lambert-app}
t_W=-L_1\big(x/\beta^\beta\big)+ \beta L_2\big(x^{1/\beta}/\beta\big)
-\beta^2\, \frac{L_2\big(x^{1/\beta}/\beta\big)}{L_1\big(x/\beta^\beta)},
\end{equation}
and by $t_C$ the approximation deduced from Comtet expansion
(\ref{t-comtet}),
\begin{equation}
\label{t-comtet-app}
t_C=-L_1(x)+ \beta L_2(x)
-\beta^2 \, \frac{L_2(x)}{L_1(x)}.
\end{equation}
\begin{table}[htb]
\centering
\begin{tabular}{lccccccc}
\toprule
&{$\bs x$} & $10^{-1}$  & $10^{-2}$ & $10^{-3} $ & $10^{-4} $  & $10^{-5} $ & $10^{-6} $ \\
\midrule
$\bs{\beta=0.5}$&$\bs t$      & 2.8212 & 5.4533 & 7.9440 &10.3803& 12.7871 & 15.1753\\
                & $\bs{t_W}$  & 2.8124 & 5.4554 & 7.9464 & 10.3824 &12.7889& 15.1768 \\
                & $\bs{t_C}$  & 2.8102 & 5.4517&  7.9440 & 10.3808& 12.7877 & 15.1759\\
\midrule
$\bs{\beta=4}$&$\bs t$&12.3607 & 15.5923& 18.6005& 21.4786& 24.2699& 26.9987\\
& $\bs{t_W}$  & 11.9175 & 15.3431 & 18.4547& 21.3922& 24.2198& 26.9717 \\
& $\bs{t_C}$  & 11.4342 & 16.0199 & 19.1148 & 21.9488 & 24.6826 & 27.3597 \\
\bottomrule
\end{tabular}
\caption{Comparison of the approximations $t_W$ and $t_C$ given
in (\ref{t-lambert-app}) and (\ref{t-comtet-app}) with the numeric solution
$t$ of equation (\ref{eq:appen}).}
\label{comparison-t}
\end{table}

\subsection*{A2. Extension of Comtet expansion}

Following the notations of Subsection \ref{comtet-extension:subsec},
 the solution $U_{\gamma,D}(x)$  of the equation
\begin{equation}
\label{eq:salvi1}
t^\gamma e^{t} D\Big(\frac{1}{t}\Big)=x,
\end{equation}
such that $t\to \infty$ when $x\to \infty$,
has an asymptotic (formal) expansion (see Robin \cite{Rob88})
\begin{equation}
\label{ud}
U_{\gamma, D}(x)=L_1(x)+\sum_{n=0}^\infty\frac{Q_n(L_2(x))}{L_1^n(x)},
\end{equation}
where the polynomials
 verify
the recurrence relation
\begin{equation}
\label{rec:pn}
Q'_{n+1}=-\gamma(Q'_n-nQ_n), \ n\ge 0, \quad\text{and}\quad{Q'_0=-\gamma}.
\end{equation}
Due that this recurrence does not determine the independent term of the
 polynomials,
Salvi \cite{Sal92} considers the generating function of the independent terms
of the polynomials, $Q_0(0),Q_1(0),\dots$,
$${\cal G}(s):=\sum_{n=0}^\infty Q_0(0)s^n,$$
and he proves that it satisfies
\begin{equation}
\label{rec:p0}
{\cal G}(s)=-\gamma \log(1+s\,{\cal G}(s))-\log D\Big(\frac{s}{1+s\ {\cal G}(s)}\Big),
\end{equation}
which allows to compute the independent terms.
Joining (\ref{rec:pn}) and (\ref{rec:p0}) the polynomials $Q_k(x)$ can be deduced iteratively.
 Salvi \cite{Sal92} gives the code of a Maple program to compute recurrently
those polynomials. The first two polynomials are given in (\ref{q2}).

\subsection*{A3. The secondary branch of the inverse function of
$\bs{t e^{t} D\big(1/t\big)}$}
As we commented in Subsection \ref{lambert-extension:subsec}, the equation
\begin{equation}
\label{eq-lambert-A}
t\,e^{t} D\Big(\frac{1}{t}\Big)=x
\end{equation} has a unique solution $t<0$  for $x<0$,
 such that $t\to -\infty$ when $x\to 0^-$; write $t=W_{-1,D}(x)$.
 Its asymptotic expansion  can be
deduced from the very generals results of Robin \cite{Rob88} and
Salvi \cite{Sal92} commented in Appendix A2:
Indeed, equation (\ref{eq-lambert-A}) is equivalent to
$$\big(-t\big)^{-1}e^{-t}D^{-1}\bigg(-\frac{1}{-t}\bigg)=-\frac{1}{x},$$
where $D^{-1}(t)=1/D(t)$. Since $-\dfrac{1}{x}\to \infty$,
changing $-t$ by $u$, we  have
$$u^{-1} e^uC(u)=-1/x,$$
where $C(t)=1/D(-t).$
Hence, with the notations of Appendix A2,
the solution is
$$t=-U_{-1,C}(-1/x).$$
Thus, from (\ref{ud}), the asymptotic expansion of $W_{-1,D}(x)$ is
\begin{align}
\label{lambert-D-extension}
W_{-1,D}(x)&=-U_{-1,C}(-1/x)=-L_1(-1/x)-\sum_{n=0}^\infty \frac{R_n(L_2(-1/x))}{L_1^n(-1/x)} \notag\\
&=L_1(-x)+\sum_{n=0}^\infty (-1)^{n+1}\frac{R_n(L_2(-x))}{L_1^n(-x)},
\end{align}
where the polynomials $R_n$ satisfy
\begin{equation}
\label{ap:rec:qn}
R'_{n+1}=R'_n-nR_n, \ n\ge 0, \quad\text{and}\quad{R'_0=1},
\end{equation}
and the generating function of the independent terms of the polynomials,
${\cal H}(s):=\sum_{n=0}^\infty R_n(0)s^n$, satisfies
\begin{equation}
\label{ap:rec:q0}
{\cal H}(s)=\log(1+s\,{\cal H}(s))-\log C\Big(\frac{s}{1+s\ {\cal H}(s)}\Big)=
\log(1+s\,{\cal H}(s))+\log D\Big(-\frac{s}{1+s\ {\cal H}(s)}\Big).
\end{equation}
Again, from (\ref{ap:rec:qn}) and (\ref{ap:rec:q0}) the polynomials $R_k(x)$
 can be computed iteratively. The first two polynomials are given in
(\ref{r2}).

\subsection*{A4. Comparison of two asymptotic expansions}
In a similar way that  in Appendix A1, we are going to compare two asymptotic
expansions  for the solution of the equation
$$t^\beta e^{-t} A\bigg(\frac{1}{t}\bigg)=x,$$
where $\beta>0$ and $A(x)=\sum_{n=0}^\infty a_n x^n,$
such that $t\to \infty$ when $x\to 0^+$.
The first way to get the solution is noting that
$$t=U_{-\beta, A^{-1}}(1/x).$$
The (formal) asymptotic expansion  is (see Appendix A2)
\begin{align}
t=L_1(1/x)+\sum_{n=0}^\infty \frac{Q_n(L_2(1/x))}{L_1^n(1/x)}
=-L_1(x)-\sum_{n=0}^\infty (-1)^{n+1}\frac{Q_n(L_2(x))}{L_1^n(x)},
\end{align}
where
\begin{equation}
\label{ap:rec:pn}
Q'_{n+1}=\beta(Q'_n-nQ_n), \ n\ge 0, \quad\text{and}\quad{Q'_0=\beta},
\end{equation}
and
${\cal G}(s):=\sum_{n=0}^\infty Q_n(0)s^n$, verifies
\begin{equation}
\label{ap:rec:p0}
{\cal G}(s)=\beta \log(1+s\,{\cal G}(s))-\log A^{-1}\Big(\frac{s}{1+s\ {\cal G}(s)}\Big).
\end{equation}
The second asymptotic expansion is deduced from
$$t=-\beta W_{-1,D}\bigg(-x^{1/\beta}/\beta\bigg),$$
where
$$D(t)=A^{1/\beta}\bigg(-\frac{1}{\beta}\,t\bigg),$$
and $A^{1/\beta}(t)=\big(A(t)\big)^{1/\beta}$.
Then, from (\ref{lambert-D-extension}),
\begin{align}
t&=-\beta\bigg(L_1\big(x^{1/\beta}/\beta\big)+\sum_{n=0}^\infty (-1)^{n+1}
\frac{R_n\Big(L_2\big(x^{1/\beta}/\beta\big)\Big)}{L_1^n\big(x^{1/\beta}/\beta\big)}\bigg)\\
&=-L_1(x^\beta/\beta^\beta)-\sum_{n=0}^\infty (-\beta)^{n+1}
\frac{R_n\Big(L_2\big(x^{1/\beta}/\beta\big)\Big)}{L_1^n(x^\beta/\beta^\beta)},
\end{align}
where the polynomials $R_n$ are determinated by (\ref{ap:rec:qn})
and (\ref{ap:rec:q0}).
With the notations of Appendix A3, it follows that $Q_n(x)=\beta^{n+1} R_n(x)$: just define the polynomials $S_n(x)=\beta^{n+1}
R_n(x)$ and check that they satisfy (\ref{ap:rec:pn}), and that the corresponding
generating function of the independent terms satisfy (\ref{ap:rec:p0}).
So, as in Appendix A1, we deduce that the two
 asymptotic  expansions  are the same function applied to different points,
and the analysis of Appendix A1  can be extended to this more general context.

\section*{Acknowledgments}
The first author was   partially
 supported by grants MINECO/FEDER reference MTM2008-03437
 and Generalitat de Catalunya reference 2009-SGR410.
The second author by the European Space Agency (ESA) under the DINGPOS contract AO/1-5328/06/NL/GLC and by the Spanish Government under project TEC2011-28219
 The  third author by grants
  MINECO/FEDER reference
   MTM2009-08869 and MINECO reference  MTM2012-33937.

\end{document}